\input amstex
\documentstyle{amsppt}
\magnification=1200

\vsize19.5cm \hsize13.5cm \TagsOnRight \pageno=1
\baselineskip=15.0pt
\parskip=3pt

\def\det{\text{det}}

\def\phi{\varphi}

\topmatter

\title  Relative $K$-stability and modified K-energy  on toric manifolds
\endtitle
\author Bin $\text{Zhou}^*$ and Xiaohua $\text{Zhu}^{*,**}$
\endauthor

\thanks {\flushpar **  Partially supported by  NSF10425102 in China and the Huo Y-T Fund.}\endthanks
\address{ * Department of Mathematics, Peking University,
Beijing, 100871, China}\endaddress
 \email{xhzhu\@math.pku.edu.cn} \endemail

\abstract
 In this paper, we discuss  the  relative $K$-stability  and
  the modified  $K$-energy associated to the Calabi's extremal
metric on toric manifolds. We give a sufficient condition in the
sense of convex polytopes  associated to  toric manifolds for both the
relative $K$-stability and the properness of modified  $K$-energy.
In particular, our results hold for toric Fano manifolds with
vanishing Futaki-invariant. We also verify our results on the
toric Fano surfaces.
\endabstract
\endtopmatter

\document

\subheading {\bf 0. Introduction }

The relation  between some geometric stabilities and the existence
of K\"ahler-Einstein metrics, more general, of extremal metrics
has been recently studied extensively ([T2], [D1], [D2], [M2],
[M3], etc.). The goal is to establish a necessary and sufficient
condition for the existence of extremal metrics in the sense of
geometric invariants ([Ya], [T2]). More recently, K.S. Donaldson
studied a relation between $K$-stability of Tian-Donaldson and
Mabuchi's $K$-energy  on a polarized toric manifold, and for the
surface's case he proved that $K$-energy is bounded from below
under the assumption of  $K$-stability ([D2]). So it is a natural
problem how to verify $K$-stability on a polarized K\"ahler
manifold $M$ although we do not know if there exists an extremal
metric on $M$. In this paper, we focus our attention on polarized
toric manifolds  as in [D2] and give a sufficient condition  for
$K$-stability in the sense of integral  polytopes associated to
polarized toric manifolds. In fact, we can study the relative
$K$-stability ( or the relative $K$-semistability) involved an
extremal holomorphic vector field related to an extremal metric on
$M$ ([Sz]).  This relative $K$-stability  is a generalization of
$K$-stability  and removes the obstruction arising from the
non-vanishing of Futaki invariant ([FM]).

An n-dimensional polarized toric manifold $M$ is corresponding to
an integral   polytope $P$ in $\Bbb R^n$ which is described by a
common set of some half-spaces,
$$\langle l_i, x\rangle < \lambda_i, ~i=1,...,d. \tag 0.1$$
 Here $l_i$ are $d$ vectors in $\Bbb R^n$ with all components in $\Bbb
 Z$, which satisfy the Delzant condition ([Ab]). Without the loss of generality,
we may assume that the original point $0$ lies in $P$, so all
$\lambda_i>0$. Let $\theta_X$ be a potential function associated
to an extremal holomorphic vector field $X$ and a K\"ahler metric
on $M$.  Define
 $$\|\theta_X\|=\max_M\{|\theta_X|\}.$$
 One can  show  that $\|\theta_X\|$ is independent of choice of K\"ahler
metrics in the  K\"ahler class ([M1]). Let $\bar R$ be the average
of scalar curvature of a K\"ahler metric in the K\"ahler class.
Then we have

\proclaim{Theorem 0.1} Let  $M$ be an $n$-dimensional polarized
toric manifold and $X$ be  an extremal vector field on $M$
associated to the polarized  K\"ahler class.  Suppose that
$$\bar{R}+\|\theta_X\| \le \frac{n+1}{\lambda_i},
~\forall ~i=1,..., d.\tag 0.2$$
 Then $M$ is $K$-stable relative to
the $\Bbb C^*$-action  induced by $X$ for any toric degeneration.
\endproclaim

According to [D2], a toric degeneration is a  special test
configuration induced by a  positive rational, piecewise  linear
function on an integral polytope associated to a polarized toric
manifold. We will show that the relative $K$-semistability for   a
toric degeneration is a necessary condition for the existence of
extremal metrics on such a manifold.  In case that $M$ is toric
Fano manifold,  condition (0.2) in Theorem 0.1 becomes
 $$\bar{R}+\|\theta_X\| \le n+1, \tag 0.2'$$
 since in this case all $\lambda_i$ are $1$.
 In particular,  (0.2') is trivial for toric Fano manifolds with vanishing Futaki-invariant since
$\theta_X=0$,  and so Theorem 0.1 is true for these manifolds. We
also verify that (0.2') is true on a toric Fano surface (cf.
Section 5). It is interesting to ask whether (0.2') is true for
higher dimensional toric Fano manifolds or not.

 The idea of proof of Theorem 0.1 is to estimate a modified Futaki-invariant
related to an extremal metric in the sense of a linear functional
on a class of convex functions on $P$. The late will rely  on the
structure of an integral polytope associated to a polarized toric
manifold (cf. Section 3, 4). Our method here can  be also used to
study a modified  $K$-energy $\mu(\phi)$ on a toric K\"ahler
manifold, not only a polarized toric manifold. In the other words,
the numbers $\lambda_i$ in the relation (0.1)  need not to be
integral. Note that the complex structure on a toric manifold $M$
is uniquely determined by the set $\{l_i\}_{i=1,...,d}$ of $d$
vectors in $\Bbb R^n$ and a different  set
$\{\lambda_i\}_{i=1,...,d}$ of $d$ numbers determines a different
K\"ahler class ([Gui]).

 As same as the $K$-energy, $\mu(\phi)$ is defined on a K\"ahler class whose
critical point is an extremal metric while the critical point of
$K$-energy is a K\"ahler metric with constant scalar curvature.

\proclaim{Definition 0.1} Let  $\omega_g$ be a K\"ahler form of
K\"ahler metric $g$ on a compact K\"ahler manifold $M$ and
$\omega_{\phi}$ be a K\"ahler form associated to  a potential
function $\phi$ in K\"ahler class $[\omega_g]$.  Let
$$I(\phi)=\frac{1}{V}\int_M \phi(\omega_g^n-\omega_{\phi}^n),$$
where $V=\int_M\omega_g^n$.   $\mu(\phi)$ is called proper
associated to a subgroup $G$ of the automorphisms group
$\text{Aut}(M)$ in K\"ahler class $[\omega_g]$  if there is a
continuous function $p(t)$ in $\Bbb R$ with the property
 $$ \lim_{t\to+\infty} p(t)=+\infty,$$
  such that
  $$\mu(\phi)\ge c\inf_{\sigma\in G} p(I(\phi_{\sigma}))-C$$
   for some uniform positive constants $c$ and $C$ independent of $\phi$, where
 $\phi_{\sigma}$ is defined by
$$\omega_g+\sqrt{-1}\partial\bar{\partial}\phi_{\sigma}
=\sigma^*(\omega_g+\sqrt{-1}\partial\bar{\partial}\phi).$$
\endproclaim

In  the following  we give  a sufficient condition   for the
properness of $\mu(\phi)$ on a toric K\"ahler manifold.

 \proclaim{Theorem 0.2} Let $M$ be a  toric K\"ahler  manifold associated to a convex polytope $P$.
Let $G_0$ be a maximal compact subgroup of torus actions group $T$
on $M$. Suppose that condition (0.2) is satisfied. Then
$\mu(\phi)$ is proper associated to subgroup $T$ in the space of
$G_0$-invariant K\"ahler metrics. In fact, there exists numbers
$\delta
>0$ and $C$ such that for any $G_0$-invariant potential functions
associated to the K\"ahler class, it holds
$$\mu(\phi) \geq \delta \inf_{\tau\in {T}}{I(\phi_\tau)} - C_{\delta}.\tag 0.3$$
\endproclaim

Since condition (0.2)  in Theorem 0.1 is true on a toric Fano
surface, we see that inequality (0.3) holds on such a  manifold.
It is known that toric Fano surfaces are classified into five
different types, i.e., $\Bbb CP^2$, $\Bbb CP^1\times \Bbb CP^1$
and $\Bbb CP^2\# l\overline {\Bbb CP^2}~(l=1,2,3)$ and there
exists a K\"ahler-Einstein on each one except $\Bbb CP^2\#
l\overline {\Bbb CP^2}~(l=1,2)$ ([T1]). On the other hand, it is
still unknown whether there exists an extremal metric on $\Bbb
CP^2\# 2\overline {\Bbb CP^2}$ while  the extremal metric was
constructed on $\Bbb CP^2\# 1\overline {\Bbb CP^2}$ by using the
ODE method ([Ca], [Ab]). We hope that our result about the
properness of modified $K$-energy  will  be useful  in the study
of the existence problem of extremal metrics
 on $\Bbb CP^2\# 2\overline {\Bbb CP^2}$ ([LS]). In fact, one may
 guess that there  always exists   an extremal metric on the K\"ahler class
 $2\pi c_1(M)$ on a toric Fano manifold.  We note that the
 existence problem was solved recently in [WZ] when $\theta_X=0$,
  which is equal to the existence of K\"ahler-Einstein metrics
 with positive scalar curvature.

When the Futaki invariant vanishes on the  K\"ahler class, the
extremal vector holomorphic  field should be  zero. Then condition
(0.2) becomes
 $$\bar{R}\le \frac{n+1}{\lambda_i},
~\text{for all},~i=1,..., d.\tag 0.2''$$
 Condition (0.2'') seems not so hard to satisfy. For instants, we construct a kind of
 of K\"ahler classes  on  the
toric manifold $\Bbb CP^2\# 3\overline{\Bbb CP^2}$, which satisfy
the condition (0.2'') (Section 6).

 Since $\bar{R}$ can be computed as some integrals in the sense of
 convex polytope $P$ (cf. Section 4),   as a corollary of Theorem
0.2, we have

\proclaim {Corollary 0.1} Let  $M$ be an $n$-dimensional  toric
K\"ahler  manifold with vanishing Futaki invariant. Suppose
 $$ \max_{i=1,...,d} \{ \frac{\lambda_i}{ \sum_{j=1}^d
\int_{P_j}dx } \sum_{j=1}^d \frac{1}{\lambda_j}\int_{P_j}dx\} \le
\frac {n+1}{n}, \tag 0.4$$
  where $\{P_j\}_{j=1}^d$ is a  union  of cones with  $(n-1)$-dimensional
 faces $E_j$ on $\partial P$ as
 bases and vertex at $0$.  Then  $K$-energy $\mu(\phi)$ satisfies (0.3). In
particular, $K$-energy $\mu(\phi)$ is bounded from below in the
space of $G_0$-invariant potential functions.
\endproclaim

Finally, we remark that  the properness of $K$-energy is a
necessary and sufficient condition for the existence of
K\"ahler-Einstein metrics with positive scalar curvature,  and in
case that $M$ is a K\"ahler-Einstein manifold with positive first
Chern class, (0.3) implies
$$\mu(\phi) \geq \delta
I(\phi) - C_{\delta},\tag 0.3'$$
  for any $\phi\in\Lambda_1^{\bot}(M,g_{KE}),$
 where  $\Lambda_1(M,g_{KE})$ denotes the first non-zero eigenfunctions space
for  the Lapalace operator associated to a K\"ahler-Einstein
metric $g_{KE}$ (cf. Section 4).  The existence of such an
 inequality was conjectured by Tian for any K\"ahler-Einstein
manifold with positive first Chern class ([T2]).

The organization of this paper is as follow. In Section 1, we
derive a Donaldson's version for the modified $K$-energy on toric
manifolds. Then in Section 2,  we give an analytic criterion for
the properness of modified $K$-energy. In Section 3, we recall the
relative $K$-stability (relative $K$-semistability) introduced in
[Sz] and compute the modified Futaki invariant for a kind of toric
degenerations.  Both Main Theorem 0.1 and 0.2  will be proved in
Section 4. In Section 5, we estimate the potential function of an
extremal vector field on $\Bbb CP^2\#\overline{ 1 \Bbb CP^2}$ and
$\Bbb CP^2\#\overline{ 2 \Bbb CP^2}$ respectively, then verify
condition (0.2) for toric Fano surfaces.  In Section 6,  we
construct a kind of  of K\"ahler classes  on  the toric manifold
$\Bbb CP^2\# 3\overline{\Bbb CP^2}$, which satisfy the condition
(0.2'').

\subheading {\bf 1. A  Donaldson's version of $\mu(\phi)$}
\vskip3mm

In this section, we review some basic materials of toric
differential geometry, then apply them to  the modified $K$-energy
$\mu(\phi)$.

Let $(M, g)$ be a compact K\"ahler manifold with dimension $n$.
Then  K\"ahler form $\omega_g$ of $g$ is given  by
$$\omega_g=\sqrt{-1}\sum_{i,j=1}^n g_{i\overline j} dz^i\wedge d\overline
z^j$$
 in local coordinates $(z_1,...,z_n)$, where $g_{i\overline j}$ are components of metric $g$.
 The K\"ahler class $[\omega_g]$ of $\omega_g$ can be represented by a set of potential
 functions as follow
 $$\Cal M =\{\phi\in C^{\infty}(M)|
~\omega_{\phi}=\omega_g+ \frac {\sqrt {-1}}{2\pi}
\partial\overline\partial\phi >0\}.$$
According to [Ca], a K\"ahler metric $\omega_{\phi}$ in the
K\"ahler class $[\omega_g]$ is called extremal if
  $$R(\omega_{\phi})=\overline R+\theta_X(\omega_\phi)$$
   for some Hamiltonian holomorphic vector field $X$ on $M$, where
   $R(\omega_{\phi})$ is the
scalar curvature of $\omega_{\phi}$,  $\overline R=\frac
{1}{V}\int_M R(\omega_g)\omega_g^n$, $V=\int_M \omega_g^n$ and
$\theta_X(\omega_\phi)$ denotes the potential function of $X$
associated to the metric $\omega_\phi$,  which is defined by
$$\cases &i_X\omega_{\phi} =  \sqrt{-1}\overline\partial \theta_X(\omega_{\phi}),\\
&\int_M \theta_X(\omega_{\phi})\omega_\phi^n=0.\endcases\tag 1.1
$$
By [FM], such an  $X$, usually called extremal is uniquely
determined by the Futaki invariant $F(\cdot)$. When the Futaki
invariant vanishes, then $X=0$ and  an extremal metric becomes one
with constant scalar curvature. It is well-known that a K\"ahler
metric with constant scalar curvature is a critical point of the
Mabuchi's $K$-energy. In the case of extremal metrics, one can
modify the $K$-energy to be
$$\mu (\phi) = - \frac{1}{V}\int_{0}^{1}\int_{M}
\dot\psi_{t}[R(\omega_{\psi_{t}}) -\overline R- \theta_X(\psi_t)]
\omega_{\psi_{t}}^n\wedge dt,$$
  where $\psi_{t} (0 \leq t \leq 1)$
is a path connecting $0$ to $\phi$ in $\Cal M$.  It can be showed
that the functional $\mu(\phi)$ is well-defined, i.e., it is
independent of the choice of path $\psi_t$ ([Gua]). Thus $\phi$ is
a critical point of $\mu(\cdot)$ iff the corresponding metric
$\omega_\phi$ is extremal.

Now we assume that $M$ is  an n-dimensional toric K\"ahler
manifold and $g$  is a $G_0\cong (S^1)^n$-invariant K\"ahler
metric in the K\"ahler class,  where $G_0$ is a maximal compact
subgroup of torus actions group $T$ on $M$. Then under an affine
logarithm coordinates system $(w_1,...,w_n)$, its K\"ahler form
$\omega_g$ is determined by a convex function $\psi_0$ on $\Bbb
R^n$, namely
$$\omega_g = \sqrt{-1}\partial\bar{\partial}\psi_0$$
is defined on the open dense orbit $T$. Denote $D\psi_0$ to be a
 gradient map (moment map) associated to $T$.
 Then the image of $D\psi_0$ is  a  convex polytope
$P$ in $\Bbb R^n$.
 By using the Legendre transformation $y=(D\psi_0)^{-1}(x)$, we see
that the function (Legendre function) defined by
$$u_0(x)=\langle y,D\psi_0(y) \rangle -\psi_0(y)=\langle y(x),x\rangle-\psi_0(y(x))$$
is  convex on $P$.  In general, for any $G_0$-invariant potential
function $\phi$  associated to the K\"ahler class $[\omega_g]$,
one  gets a convex function $u(x)$ on $P$ by using the above
relation  while $\psi_0$ is replaced by $\psi_0+\phi$. Set
 $$\Cal C =\{u=u_0+v| ~ u~ \text{is a convex function in }~ P,~v\in
C^{\infty}(\bar{P})\}.$$
 It was showed in [Ab] that  functions in $\Cal C$ are
 corresponding to
$G_0$-invariant functions in $\Cal M$ (whose set is  denoted by
$\Cal M_{G_0}$)  by one-to-one.

The convex polytope $P$ is described by a common set of some half-spaces,
 $$\langle l_i, x\rangle < \lambda_i, ~i=1,...,d,\tag 1.2$$
 where $l_i$ are d-vectors in $\Bbb R^n$ with all components in $\Bbb
 Z$, which satisfy the Delzant condition. Conversely, given a convex polytope $P$
in $\Bbb R^n$ as above, one can construct a $G_0$-invariant
K\"ahler metric on $M$ ([Gui]). Without the loss of generality, we
may assume that the original point $0$ lies in $P$, so all
$\lambda_i>0$. Let $d\sigma_0$ be the Lebesgue measure on the
boundary $\partial P$ and $\nu$ be the outer normal vector field
on $\partial P$. Let $d\sigma=\lambda_i^{-1}(\nu, x)d\sigma_0$ on
the face $\langle l_i, x\rangle= \lambda_i$ of $P$. It is clear
that $(\nu, x)$ is constant on each face. According to [D2], if $u
\in \Cal C$ corresponds to $\phi \in \Cal M_{G_0}$, then
$$ -\frac{1}{2^nn!}\int_{0}^{1}\int_{M}
\dot\psi_{t}R(\omega_{\psi_{t}}) \omega_{\psi_{t}}^n
dt=(2\pi)^n\left(-\int_{P}\log(\det{D^2u})dx + \int_{\partial P}u
d\sigma\right).$$
 By using a  relation
$$\dot{\psi_{t}} = - \dot{u_{\psi_t}},$$
 we get
  $$\mu(\phi)=\frac{2^nn!(2\pi)^n}{V} \Cal F(u),$$
where
$$\Cal F(u) = - \int_{P}\log(\det{D^2u})dx + \int_{\partial P}u d\sigma
- \int_{P}(\bar{R}+\theta_X)udx.\tag 1.3$$

\proclaim {Lemma 1.1}
 $$\theta_X=\sum_{i=1}^n a_i(x_i + c_i),$$
  where $2n$-constants $a_i$ and $c_i$ are determined uniquely by $2n$-equations,
$$\frac{\text{vol}(P)}{\text{vol}(M)}F(\frac{\partial}{\partial y_i})=- \int_P (\sum_{j=1}^n a_j(x_j + c_j))(x_i+c_i)
dx,~i=1,...,n,\tag 1.4$$
$$\int_P (x_i+c_i) dx=0,~i=1,...,n.\tag 1.5$$
\endproclaim

\demo{Proof} According to [FM], an extremal holomorphic field $X$
should be belonged to the center of a  reductive Lie subalgebra of
Lie algebra consisting of holomorphic vector fields on $M$. In
particular, $X$ is an element of an Abelian Lie subalgebra. So
under the affine coordinates system, we see
$$X\in\text
{span}\{\frac{\partial}{\partial y_i},~i=1,...,n\}.$$
 Note that
$$\theta_{\frac{\partial}{\partial y_i}}=x_i+c_i\tag 1.6$$
 for each $\frac{\partial}{\partial y_i}$, where  constants $c_i$ are uniquely
 determined by (1.5) according to (1.1).
Thus
 $$\theta_X=\sum_{i=1}^n a_i(x_i + c_i)\tag 1.7$$
 for some constants $a_i$. Since  $\theta_X$ satisfies ([FM]),
$$\frac{\text{vol}(P)}{\text{vol}(M)} F(\frac{\partial}{\partial y_i})=-\int_P
\theta_X\theta_{\frac{\partial}{\partial y_i}}dx, ~i=1,...,n,\tag
1.8$$
 then by relation (1.6) and (1.7), we get (1.4).\qed \enddemo

\proclaim {Lemma 1.2} Let $u$ be an affine linear function on $P$.
Then
$$L(u)= \int_{\partial P}u d\sigma -
\int_{P}(\bar{R} +\theta_X)u dx=0.$$
  \endproclaim

\demo {Proof}  First note ([D2])
 $$\bar R=\frac{\int_P dx}{\int_{\partial P} d\sigma}.\tag 1.9$$
 Then for any constant $c$, we have $L(c)=0.$ Thus we suffice to
prove
 $$ L(x_i)=0, ~i=1,...,n.$$
Let $\rho_t^i$ be a one-parameter subgroup on $M$ generated by
$\text{re}(\frac{\partial}{\partial y_i})$. Then
$$ \rho_t^i(y)= \cases &y_j,~j\neq i\\
&y_j+t,~j=i.\endcases$$
 Thus $(\phi+\psi_0)_{\rho_t^i}(y)=(\phi+\psi_0)(\rho_t^i(y))$ and
$$u_{\rho_t^i}(x)=u- tx_i,$$
 where $(\phi+\psi_0)_{\rho_t^i}$ is a convex function in $\Bbb R^n$
induced by actions $\rho_t^i$ as follow,
$$\sqrt{-1}\partial\overline\partial(\phi+\psi_0)_{\rho_t^i}=(\rho_t^i)^*\omega_{\phi},$$
 and $u_{\rho_t^i}$ is the Legendre function of
$(\phi+\psi_0)_{\rho_t^i}$.
 Since
$$ \text{vol}(M)\frac{\mu(\phi_{\rho_t^i})}{dt}|_{t=0}=F(\frac{\partial}{\partial y_i})
+\int_M\theta_X \theta_{\frac{\partial}{\partial
y_i}}\omega_g^n=0,$$
 we get
$$ 2^n n!(2\pi)^n \frac{\Cal F(u_{\rho_t^i})}{dt}|_{t=0}=-L(x_i)=0.$$
 \qed
\enddemo

 By Lemma 1.2,  we see that the
 functional $\Cal F(u)$ is invariant if $u$ is replaced by adding
 an affine linear function.  For this reason, we normalize $u$ as
follows. Let  $p\in P$ and set
 $$\tilde{\Cal C}= \{u \in \Cal C
|\inf_{P}u=u(p)=0 \}.$$
 Then for any $u_{\phi} \in \Cal C$ corresponding to $\phi\in \Cal M_{G_0}$, one can normalize $u_{\phi}$ by
$$\tilde{u}_{\phi} = u_{\phi} - (\langle Du_{\phi}(p),x-p\rangle + u_{\phi}(p))$$
so that $\tilde{u}_{\phi}=u_{\tilde\phi}\in \tilde{\Cal C}$
corresponds to a K\"ahler potential function $\tilde{\phi}\in \Cal
M_{G_0}$  which satisfies
$$D(\tilde{\phi} + \psi_0)(0) = p~\text{and}~ (\tilde\phi+\psi_0)(0)=0.\tag 1.10$$
In fact, $\tilde\phi$ can be uniquely determined by using the
affine coordinates transformation $y\to y+y_0$ as follow,
$$\tilde\phi(y)=(\phi+\psi_0)(y+y_0)-\psi_0(y)- (\phi+\psi_0)(y_0).$$

\subheading { \bf 2 An analytic criterion for the properness of
$\mu(\phi)$}
 \vskip3mm

In this section, we give an analytic criterion for the properness
of the modified $K$-energy $\mu(\phi)$.
 We need to recall the Aubin's functional ([Au]),
$$J(\phi) =
\frac{1}{V}\int_{0}^{1} \int_{M}
\dot{\psi_{t}}(\omega^{n}_g-\omega^{n}_{\psi_t})\wedge dt,
~\forall ~\phi\in \Cal M,
$$
 where  $\psi_t (0\le t\le 1)$ is a path  in $\Cal M$ connecting $0$ to
 $\phi$. Note
$\dot{\psi_{t}} = - \dot{u_{t}}$, where $u_t$ are Legendre
functions associated to potential functions $\psi_t$. Then
$$J(\phi)= \frac{1}{V}\int_{M} \phi \omega^{n}_g +
H(u_{\phi})-H(u_{0}),\tag 2.1$$
 where
  $$H(u) = \frac{1}{Vol(P)}\int_{P}u dx,~\forall~u\in \Cal C.$$

\proclaim{Lemma 2.1} There exists $C>0$ such that
$$|J(\tilde{\phi})- H(u_{\tilde\phi})|\leq C,
 ~\forall~ \phi\in \Cal M_{G_0},$$
where $\tilde\phi$ is as defined in (1.10).
\endproclaim

\demo{Proof}By (2.1), we have
$$J(\tilde{\phi})- H(u_{\tilde{\phi}}) =
\frac{1}{V}\int_{M}\tilde{\phi} \omega^{n}_g-H(u_{0}).$$
 We claim that
 $$|\frac{1}{V}\int_{M}\tilde{\phi} \omega^{n}_g| \le C$$
  for some uniform constant $C$.

  Applying the Green's formula to the potential function $\tilde\phi$,
    one sees that there exists a constant $C_{0}$ such that
$$\frac{1}{V}\int_{M}\tilde{\phi}\omega^{n}_g
\geq \sup{\{\tilde{\phi}\}} - C_{0}.\tag 2.2$$
  Set
  $$\Omega_N =\{  \xi\in M|~
 \tilde{\phi}(\xi)\leq \sup_{\Bbb R^n}{\{\tilde{\phi}\}} - N
\}.$$
 Then
$$ \aligned &
\frac{1}{V}\int_{M}\tilde{\phi}\omega^{n}_g\\
& = \frac{1}{V} \int_{M \cap \Omega_N}\tilde{\phi}\omega^{n}_g+
\frac{1}{V}\int_{M
\setminus \Omega_N}\tilde{\phi}\omega^{n}_g\\
& \leq  \frac{1}{V}[(\sup_{\Bbb R^n}{\{\tilde{\phi}\}} - N) Vol(M
\cap \Omega_N) + \sup_{\Bbb R^n}{\{\tilde{\phi}\}}Vol(M\setminus \Omega_N)]\\
& = \sup_{\Bbb R^n}{\{\tilde{\phi}\}} - \frac{N Vol(M \cap
\Omega_N)}{Vol(M)}.
\endaligned$$
It follows
$$Vol(M \cap \Omega_N) \leq \frac{C_{0}Vol(M)}{N}=\frac{C_0V}{N}\to 0,\tag 2.3$$
 as $N\to \infty$.
On  the other hand, by the second relation in (1.10), we have
  $$\tilde\phi(0)=-\psi_0(0).$$
Then
 $$\tilde \phi(x)\le \phi(0)-2r\sup\{|p|: ~p\in P\}\le C(r),~\forall~x\in B_r(0),$$
  where $C(r)$ depends only on the radius $r$ of ball $B_r(0)$ centered at the original.
  Since  the volume of domain $B_1(0)\times (2\pi)^n$ associated the metric
  $\omega_g=\sqrt{-1}\partial\overline\partial \psi_0$ is bigger than some number $\epsilon>0$,
  by (2.3), it is easy to see that there is at least a point $x_0 \in B_1(0)$ such that
  $$\tilde\phi(x_0)\ge \sup_{\Bbb R^n}\tilde\phi-N$$
 as $N$ is sufficiently large. Thus
  $$ \sup_{\Bbb R^n}\tilde\phi\le N+C(1),$$
and consequently
$$\frac{1}{V}\int_{M}\tilde{\phi}\omega^{n}_g\le N+C(1).$$
 By (2.2), we also get
$$\frac{1}{V}\int_{M}\tilde{\phi}\omega^{n}_g\geq \tilde\phi(0)-C_0=-\psi(0)-C_0.$$
 Therefore the claim  is true and Lemma is proved.\qed
\enddemo

As in [D2], one can extend  the definition space $\Cal C$ of
functional $\Cal F(u)$ to
$$\Cal C_{\infty}=\{u\in C^{\infty}(P)\cap C(\overline
P)|~u~\text{is convex in}~ P\}$$
 and one can show that
 $$\inf_{u\in \Cal C}\Cal F(u)=\inf_{u\in \Cal C_{\infty}}\Cal
 F(u).$$
 Note that the minimal point  of $\Cal F(u)$ is unique
  by using the convexity of $-\text{log}\text{det}$ and it satisfies the
  Euler-Langrage equation,
  $$\sum_{i,j=1}^n u^{ij}_{ij}=-(\overline R+\theta_X),$$
    where $(u^{ij})=(u_{ij})^{-1}$ and $u^{ij}_{kl}$ denote the second derivatives
   of $u^{ij}$. In fact, the above statement also holds for more general functional $\Cal F(u)$
    while the function $\overline R+\theta_X$ in (1.3) is  replaced by
    a smooth one in $\overline P$.

   The following result is due to [D2].

\proclaim{Lemma 2.2} Suppose that  there exists a $ \lambda
> 0$ such that
$$L(u)  \geq \lambda \int_{\partial P}u d\sigma, \tag 2.4$$
for any normalized function   $ u\in\Cal C_{\infty}$.  Then  there
exists a $\delta>0$ depending only on $\lambda$ such that
$$\Cal F(u)\geq \delta H(u) - C_{\delta}\tag 2.5$$
for any normalized function $u\in \Cal C_{\infty}$.
\endproclaim

\demo{Proof} Choose a function $v_0$ in $\Cal C$ and define a
smooth function $A$ in $\overline P$ by
$$v^{ij}_{0ij}=-(\overline R+A).$$
 Let $\Cal F'(u)$ be a modified functional of $\Cal F(u)$
  while the part of linear functional $L(u)$ of $\Cal F(u)$ is replaced   by
   $$L'(u)=  \int_{\partial P}u d\sigma -\int_{P} (\overline R +A)u dx.$$
Then
 $$\Cal F'(u)\ge \Cal F'(v_0)=-C_0,~\forall ~u\in \Cal C_{\infty}.\tag 2.6$$

We compute the difference between  the linear parts $L'(u)$ and
$L(u)$. Pick a $\delta > 0 $.  Note that
$$\int_{P}u dy \leq C\int_{\partial P}u d\sigma$$
because of the convexity of  $u\ge 0$. Then by the condition
(2.4), we have
  $$\aligned
|L(u)-L'(u)| & = \left|\int_{P}\theta_{X} u dx - \int_{P}Au
dx\right | \\
& \leq C_1 \int_{P} udx \\
& = C_1\left[(1+ \delta)\int_{P} u dx - C_1\delta\int_{P} udx\right]\\
& \leq C'_{\delta}\int_{\partial P} u d\sigma -
C_1\delta\int_{P} u dx\\
& \leq C'_{\delta,\lambda}L(u)- C_1\delta\int_{P} u dx.
\endaligned$$
 It follows
$$(C'_{\delta,\lambda}+ 1)L(u) \geq L'(u)+ C_1\delta\int_{P} u dx.$$
 Thus  by (2.6),  we get
$$\aligned
\Cal F(u)
& = - \int_{P}\log (\det (u_{ij})) dx + L(u) \\
& \geq - \int_{P}\log (\det (u_{ij}))dx +
\frac{L'(u)}{C'_{\delta,\lambda}+ 1} +
\frac{C_1 \delta}{C'_{\delta,\lambda}+ 1}\int_{P} u dx\\
& = \Cal F'(\frac{u}{C'_{\delta,\lambda}+ 1 })+
\frac{C_1\delta}{C'_{\delta,\lambda}+ 1}\int_{P} u dx -n \text{log}(C'_{\delta,\lambda}+ 1)\\
& \geq \frac{C_1\text{Vol}(P) \delta}{C'_{\delta,\lambda}+ 1} H(u)
- C_0'.
\endaligned$$
Replacing $\frac{C_1\text{Vol}(P)\delta}{C'_{\delta,\lambda}+ 1}$
by $\delta$, we obtain (2.5).
 \qed
\enddemo

\proclaim{Proposition 2.1} Suppose that (2.4) is satisfied. Then
there exists a number $\delta>0$ such that for any $G_0$-invariant
$\phi\in \Cal M_{G_0}$ it holds
$$\mu(\phi) \geq \delta \inf_{\tau\in {T^n}}{I(\phi_\tau)} - C_{\delta},\tag 2.7$$
where
$$\omega_g+\sqrt{-1}\partial\bar{\partial}\phi_{\tau}
=\tau^*(\omega_g+\sqrt{-1}\partial\bar{\partial}\phi).\tag 2.8$$
 In particular, $\mu(\phi)$ is bounded from below in $\Cal M_{G_0}$.
\endproclaim

\demo {Proof}  Let   $\phi\in \Cal M_{G_0}$.  Then there exists a
$\sigma \in T$ such that the Legendre function $u_{\phi_{\sigma}}$
associated to $\phi_{\sigma}$ is belonged to $\tilde{\Cal C}$.
 By Lemma 2.2, we see that
$$\mu(\phi_{\sigma}) \geq \delta H(\phi_\sigma) -
C_{\delta}.$$
 Note that
 $$\mu(\phi)=\mu(\phi_{\sigma}).$$
  Thus by Lemma 2.1, we get
$$\aligned \mu(\phi)&=\mu(\phi_{\sigma}) \geq \delta J(\phi_\sigma)
- C_{\delta}'\geq \delta \inf_{\tau\in {T^n}}J(\phi_\tau) -
C_{\delta}'\\
&\ge\frac{\delta}{n} \inf_{\tau\in {T^n}}I(\phi_\tau) -
C_{\delta}'.\endaligned $$
 Here at the last inequality we used the fact ([Au]),
  $$J(\phi)\ge \frac {1}{n}I(\phi),~\forall~\phi\in\Cal M.$$
  The proposition is proved.\qed
\enddemo

\subheading {\bf 3. Computation of the modified  Futaki-invariant
} \vskip3mm

 In this section, we recall the notation of relative
$K$-stability introduced in [Sz] and compute the modified Futaki
invariant   for a test configuration in the sense of Donaldson
([D2]). We assume that $(M,L)$ is  an n-dimensional polarized
toric manifold and $g$  is a $G_0\cong (S^1)^n$-invariant K\"ahler
metric with its K\"ahler form $\omega_g\in 2\pi c_1(L)$ on $M$ ,
where $L$ is a positive holomorphic line bundle on $M$. In the
other words, the corresponding convex polytope $P$ induced by the
moment map is integral (cf. Section 2).

\proclaim{Definition 3.1 ([Do], [Sz])} A test configuration for a
polarized variety  $(M,L)$ of exponent $r$ consists of a $\Bbb
C^*-$equivariant flat family of schemes $\pi: \Cal
W\longrightarrow \Bbb C$ (where $\Bbb C^*$ acts on $\Bbb C$ by
multiplication) and a $\Bbb C^*-$equivariant ample line bundle
$\Cal L$ on $\Cal W$. We require that the fibres $(\Cal W_t, \Cal
L|_{\Cal W_t})$ are isomorphic to $(M,L^r)$ for any  $t\neq 0$. A
test configuration is called trivial if $\Cal W=M \times \Bbb C$
is a product.

If $(M, L)$ is equipped with a $\Bbb C^*-$action $\beta$, we say
that a test configuration is compatible with $\beta$, if there is
a $\Bbb C^*-$action $\tilde{\beta}$ on $(\Cal W,  \Cal L)$ such
that $\pi:\Cal W\longrightarrow \Bbb C$ is $\tilde{\beta}$
equivariant with trivial $\Bbb C^*-$action on $\Bbb C$ and the
restriction of $\tilde{\beta}$ to $(\Cal W_t, \Cal L|_{\Cal W_t})$
for nonzero $t$ coincides with that of $\beta$ on $(M, L^r)$ under
the isomorphism.
\endproclaim

Note that a $\Bbb C^*$-action on $\Cal W$  induces  a $\Bbb
C^*$-action on the central fibre $M_0=\pi^{-1}(0)$ and the
restricted line bundle $\Cal L|_{M_0}$. We denote by
$\tilde{\alpha}$ and $\tilde{\beta}$ the induced $\Bbb
C^*-$actions of $\alpha$ and $\beta$ on $(M_0,\Cal L|_{M_0})$,
respectively. The relative $K$-semistability is based  on the
following modified Futaki invariant on the central fibre,
$$F_{\tilde{\beta}}(\tilde{\alpha})= F(\tilde{\alpha})-
\frac{(\tilde{\alpha}, \tilde{\beta})}{(\tilde{\beta},
\tilde{\beta})}F(\tilde{\beta}),\tag 3.1$$
  where
$F(\tilde{\alpha})$ and $F(\tilde{\beta})$ are generalized Futaki
invariants of $\tilde{\alpha}$ and $\tilde{\beta}$ defined in
[D2], respectively, and $(\tilde{\alpha}, \tilde{\beta})$ and
$(\tilde{\beta}, \tilde{\beta})$ are inner products defined in
[Sz] (also to see (3.7) below).

\proclaim{Definition 3.2 ([Sz])} A polarized variety $(M,L)$ with
a $\Bbb C^*$-action $\beta$ is $K$-semistable relative to $\beta$
if $F_{\tilde{\beta}}(\cdot)\leq 0$ for any test-configuration
compatible with $\beta$. It is called relative $K$-stable in
addition that
 the equality holds if and only if the test-configuration is trivial.
\endproclaim

Now we consider for a  polarized toric manifold $(M, L)$ which
corresponds to an integral polytope $P$ in $\Bbb R^n$. Recall that
a piecewise linear (PL) function  $u$ on $P$ is a form of
$$u= \text{max}\{u^1,...,u^r\},$$
 where $u^\lambda=\sum a_i^\lambda x_i + c^\lambda, ~\lambda=1,...,r,$
for some vectors $(a^\lambda_i)\in \Bbb R^n$ and some numbers
$c^\lambda\in \Bbb R$. $u$ is called a rational PL-function if
components $a^\lambda_i$ and numbers $c^\lambda$ are all rational.
For a  rational PL function $u$ on $P$,  choose an integer $R$ so
that
$$Q=\{(x,t)|~x\in P, 0<t<R-u(x)\}$$
is a convex polytope $Q$ in $\Bbb R^{n+1}$.  Then $M_Q$
corresponds to
 $(n+1)$-dimensional toric variety and  $L$ on $M$ induces a holomorphic  line bundle $\Cal L$
 on $M_Q$ by using the natural embedding $i: M\rightarrow M_Q$. Decomposing
 a torus action $T_{\Bbb C}^{n+1}$ on $M_Q$ as $T_{\Bbb C}^n \times \Bbb C^* $ so that
$T_{\Bbb C}^n\times \{\text{Id}\}$ is isomorphic to the torus
action on $M$, we get a $\Bbb C^*-$action $\alpha$ by
$\{\text{Id}\}\times{\Bbb C^*}$, and so we define an equivariant
map
$$\pi: M_Q\rightarrow \Bbb CP^1$$
 satisfying $\pi^{-1}(\infty)=i(M)$.
Then $\Cal W=M_Q\backslash i(M)$ is  a test configuration for the
pair $(M,L)$, called a toric degeneration ([D2]).  This test
configuration  is compatible to an extremal $\Bbb C^*-$action
$\beta$ induced by an extremal holomorphic vector field $X$ on
$M$. In fact, the extremal $\Bbb C^*-$action $\beta$ is isomorphic
to a one parameter subgroup of $T_{\Bbb C}^n \times
\{\text{Id}\}$, which acts on $\Cal W$. Since the action is
trivial in the direction of $\alpha$, it is compatible. To compute
the modified Futaki invariant for such a test configuration, we
need

\proclaim{Lemma 3.1}  Let
 $B_{k,P}=Z^n\bigcap k\bar{P}$ for
any $k\in{Z^+}$.  Assume $P$ is an integral polytope in $R^n$ and
$\phi$, $\psi$ are  two positive rational, PL-functions on $P$.
Then
$$\sum_{I\in{B_{k,P}}}\phi(I)=k^{n}\int_P\phi dx +
\frac{k^{n-1}}{2}\int_{\partial P}\phi d\sigma +O(k^{n-2})\tag
3.2$$
  and
$$\sum_{I\in{B_{k,P}}}\phi(I)\psi(I)=k^n\int_P\phi\psi dx+O(k^{n-1}).\tag 3.3$$
\endproclaim

\demo {Proof}The relation (3.3) in the lemma is trivial and it
suffices to prove (3.2). Let $Q$ be a convex polytope associated to
rational, PL-function $\phi$ as above and $B_{k,Q}=Z^{n+1}\bigcap
k\bar{Q}$. Then it is easy to see that
$$\sum_{I\in{B_{k,P}}}\phi(I)=N(B_{k,Q})-N(B_{k,P}),$$where
$N(B_{k,Q}),N(B_{k,P})$ are the numbers of point, respectively.
 Applying Proposition 4.1.3 in [D2] for each convex polytope
P and Q, one will get (3.2).\qed\enddemo

 \proclaim{Proposition 3.1} For the above test
configuration induced by a rational PL-function $u$, we have
$$F_{\tilde{\beta}}(\tilde{\alpha})=-\frac{1}{2Vol(P)}\left(\int_{\partial P}u d\sigma-
\int_P (\bar R+\theta_X)u dx  \right).\tag 3.4$$
\endproclaim

\demo{Proof} As in [D2], we consider the space $H^0(\Cal W, \Cal
L^k)$ of holomorphic sections over $\Cal W$, which has a basis
$\{S_{I,i}\}$, where $I$ is a lattice in $B_{k,P}$ and $0\leq
i\leq k(R-u)(I)$. By using the exact sequence for large $k$,
$$0\longrightarrow H^0( \Cal W, \Cal L^k\otimes \pi^*(\vartheta(-1)))\longrightarrow H^0( \Cal W, \Cal L^k)
\longrightarrow H^0(M_0, \Cal L^k) \longrightarrow 0,$$
  $H^0(M_0, \Cal L^k)$ has a basis
$\{S_{I,k(R-u)(I)}|_{M_0}\}_{I\in{B_{k,P}}}$.
 Then by (3.2), one obtains
$$\aligned d_k&=\text{dim} H^0(M, \Cal L^k)=\text{dim} H^0(M, \Cal
L^k)=N(B_{k,P})\\
&=k^n Vol(P)+ \frac{k^{n-1}}{2}Vol(\partial P)+
O(k^{n-2}).\endaligned$$

By choosing a suitable coordinates system,  we may assume that
$$\int_P x_i dx=0.$$
Then  $\theta_X=\langle\theta, x\rangle$ for some vector $\theta$
in $\Bbb R^n$  and the one parameter subgroup $\beta$ induced by
$X$ in $T_{\Bbb C}^{n+1}$ is a  form of
$$(e^{\theta_1z}, ..., e^{\theta_nz},1), ~z\in \Bbb C^*, $$
 which act on
$S_{I,i}$ with weight $k\langle\theta,I\rangle$. On the other
hand, $\alpha$  in $T_{\Bbb C}^{n+1}$ is a form of
$$(1, ..., 1, e^z),$$
  which act on $S_{I,i}$ with weight
$k(R-u)(I)$. Thus the infinitesimal generators $A_k$ and $B_k$ of
$\tilde{\alpha}$ and $\tilde{\beta}$ are $ (d_k \times d_k)$
diagonal matrices
$$diag(... ,k(R-u)(I), ...)$$
and
 $$diag(..., k\langle I, \theta\rangle ,...),$$
 respectively.    By Lemma 3.1,  we get follows,
$$\aligned
Tr(A_k)&=\sum_{I\in{B_{k,P}}}k(R-u)(I)\\
&=k^{n+1}\int_P (R-u) dx+ \frac{k^n}{2}\int_{\partial P}(R-u)
d\sigma +O(k^{n-1}),\endaligned$$
$$\aligned
Tr(B_k)&=\sum_{I\in{B_{k,P}}}k\langle\theta,I\rangle\\
&=k^{n+1}\int_P \theta_X dx+ \frac{k^n}{2}\int_{\partial
P}\theta_X d\sigma +O(k^{n-1}),
\endaligned$$
$$\aligned
Tr(A_kB_k)&=\sum_{I\in{B_{k,P}}}k^2(R-u)(I)\langle\theta,I\rangle\\
&=k^{n+2}\int_P (R-u)\theta_X dx + O(k^{n+1}),
\endaligned$$
 and
$$\aligned
Tr(B_k^2)&=\sum_{I\in{B_{k,P}}}k^2\langle\theta,I\rangle^2\\
&=k^{n+2}\int_P \theta_X^2 dx + O(k^{n+1}).
\endaligned$$
Hence
$$ F(\tilde{\alpha})=-\frac{1}{2Vol(P)}\left[\int_{\partial P}u d\sigma
-\frac{Vol(\partial P)}{Vol(P)}\int_P u dx\right],$$
$$F(\tilde{\beta})=-\frac{1}{2Vol(P)}\int_P \theta_X^2 dx. \tag 3.5$$
Recall that
 $$\overline R=\frac{Vol(\partial P)}{Vol(P)}. $$
Therefore, we get
$$F(\tilde{\alpha})=-\frac{1}{2Vol(P)}\left[\int_{\partial P}u d\sigma
-\bar R\int_P u dx\right].\tag 3.6$$

By Lemma 1.2, we have
$$L(\theta_X)=0,$$
 which is equal to
    $$\int_{\partial P}\theta_X d\sigma -
 \int_{P}(\bar R+\theta_X )\theta_X dx=0.$$
Note that
$$\int_P \theta_X dx =0.$$
It follows
$$\int_{\partial P}\theta_X d\sigma =\int_{ P}\theta_X^2 dx.$$
 Thus by the relation ([Sz]),
$$Tr(A_kB_k)-\frac{Tr(A_k)Tr(B_k)}{d_k}
=(\tilde{\alpha}, \tilde{\beta})k^{n+2}+O(k^{n+1}),\tag 3.7$$
 we get
$$(\tilde{\alpha},\tilde{\beta})=-\int_P
\theta_X u dx. \tag 3.8$$
 Similarly, we have
$$(\tilde{\beta},\tilde{\beta})=-\int_P \theta_X^2 dx. \tag 3.9$$
Substituting (3.5), (3.6), (3.8) and (3.9) into (3.1),  we obtain
(3.4) \qed
\enddemo

\subheading {\bf 4. Proof of Theorem 0.1 and 0.2}
  \vskip3mm

In this section, we prove the main theorems in Introduction. As in
Section 1, we let  $P$ be   a convex polytope in $\Bbb R^n$
defined by (1.2) which satisfies Delzant condition,  and $L(u)$ be
the linear functional
$$L(u)= \int_{\partial P}u d\sigma -
\int_{P}(\bar{R} +\theta_X)u dx$$
 defined on the space $\Cal C$. Without the loss of generality, we may assume that
$P$ contains the original point $0$.   Let $\{E_i\}_{i=1}^d$ be a
union of $(n-1)$-dimensional faces on $\partial P$ and
$\{P_i\}_{i=1}^d$ be a union of cones with bases $E_i$ and the
vertex  at $0$. First we observe

\proclaim{Lemma 4.1}
  $$ L(u)=\sum_{i=1}^d\int_{P_i}\left[
\frac{\sum_{j=1}^n{x_{j}u_{j}}}{\lambda_i}
+(\frac{n}{\lambda_i}-\bar{R}-\theta_X)u\right]dx.\tag 4.1$$
\endproclaim

\demo{Proof} Recall  that
$$d\sigma=\frac{(\nu, x)}{\lambda_i}d\sigma_0,~\text{ on each}~ E_i,$$
 where $\nu$ is the outer normal vector field on $\partial P$.  Since
$$(\nu, x)\equiv 0, ~\text{ on each}~ P_i\setminus E_i,$$
  by the Stoke's formula, we have
$$\aligned
\int_{E_i}u d\sigma&=\int_{\partial P_i}u\frac{(\nu, x)}{\lambda_i}d\sigma_0\\
&=\int_{P_i}div(\frac{xu}{\lambda_i})dx\\
&=\int_{P_i}
(\frac{\sum_{j=1}^n{x_{j}u_{j}}}{\lambda_i}+\frac{n}{\lambda_i}u)dx.
\endaligned$$
Summing these identities, we get
 $$  \int_{\partial P}u d\sigma =\sum_{i=1}^d
\int_{P_i}(\frac{\sum_{j=1}^n{x_{j}u_{j}}}{\lambda_i}+\frac{n}{\lambda_i}u)dx.$$
Then (4.1) follows from the above.
  \qed
\enddemo

\proclaim {Lemma 4.2} Let $u$  be a normalized function at the
original point $0$. Then
 $$L(u)\ge
\sum_{i=1}^d\int_{P_i}(\frac{n+1}{\lambda_i}-\bar{R}-\theta_X)u
dx.\tag 4.2$$
\endproclaim

\demo{Proof} By Lemma 4.1, we have
$$\aligned L(u)&=\sum_{i=1}^d\int_{P_i}
\frac{\sum_{j=1}^n{x_{j}u_{j}}-u}{\lambda_i} dx\\
&+\sum_{i=1}^d\int_{P_i} (\frac{n+1}{\lambda_i}-\bar{R}-\theta_X)u
dx.\endaligned$$
 Note that  $(\sum_{i=1}^n {x_{i}u_{i}}-u)$ is the Legendre function of $u$
 and so it is nonnegative by the
the normalized condition. Thus (4.2) is true. \qed\enddemo

Note that the functional $L(u)$ can be defined for a  PL-function
$u$ on $P$ introduced in Section 3 and  (4.1) still holds.  Now we
begin to prove Theorem 0.1 and need

\proclaim{Lemma  4.3} Let  $M$ be a   toric K\"ahler manifold
associated to a convex polytope $P$ and $X$ be its corresponding extremal
vector field on $M$. Suppose that for each $i=1,...,d$, it holds
$$\overline R+\theta_X\leq \frac{n+1}{\lambda_i}, ~\text{in}~
P.\tag 4.3$$
  Then for any PL-function $u$ on $P$, we have
$$L(u)\geq 0.\tag 4.4$$
 Moreover the equality  holds if and only if $u$ is  an affine linear function.
\endproclaim

\demo{Proof}  Let $u$ be a  form of
 $$u=\text{max}\{u^1,...,u^r\},$$
  where $u^\lambda=\sum a_i^\lambda x_i + c^\lambda, ~\lambda=1,...,r,$
for some   vectors $(a^\lambda_i)\in \Bbb R^n$ and  some numbers
$c^\lambda\in \Bbb R$.   By adding
 an suitable affine linear function so that $u$ is normalized
  to be $\tilde u$  with properties
  $$\tilde u\ge 0,~\text{in}~ P$$
  and
  $$\tilde u(0)=0.$$
   Denote $\tilde u^{\lambda}$ and $\tilde c^\lambda$ to be the corresponding linear functions and numbers
   respectively, as to $u$.  Then it is easy to see that $\tilde c^\lambda\leq 0$
for all $\lambda$. Dividing  $P$  into $r$ pieces $P^1, ..., P^r$
so that each $\tilde u=\tilde u^\lambda$  is defined on each
$P^\lambda$, then by Lemma 4.1 and the condition (4.3), we have
$$\aligned
L(u)&=L(\tilde u)=\sum_{i=1}^d\int_{P_i}\left[
\frac{\sum_{j=1}^n{x_{j}\tilde{u}_{j}}}{\lambda_i}
+(\frac{n}{\lambda_i}-\bar{R}-\theta_X)\tilde{u}\right]dx\\
&=\sum_{i=1}^d\sum_{\lambda=1}^r\int_{P_i\cap P^\lambda}\left[
\frac{\sum_{j=1}^n{x_{j}\tilde{u}_{j}}}{\lambda_i}
+(\frac{n}{\lambda_i}-\bar{R}-\theta_X)\tilde{u}\right]dx\\
&=\sum_{i=1}^d\sum_{\lambda=1}^r\int_{P_i\cap P^\lambda}\left[
\frac{-\tilde c^\lambda}{\lambda_i}
+(\frac{n+1}{\lambda_i}-\bar{R}-\theta_X)\tilde u^\lambda\right]dx\\
&=\sum_{i=1}^d\sum_{\lambda=1}^r\int_{P_i\cap P^\lambda}
(\frac{n+1}{\lambda_i}-\bar{R}-\theta_X)\tilde u^\lambda\geq 0.
\endaligned\tag 4.5$$
Thus  we prove  (4.4). Moreover, since each set
$\{\frac{n+1}{\lambda_i}-\bar{R}-\theta_X=0\}\cap P$ lies  in a
hyperplane in $\Bbb R^n$, one sees that   the equality in (4.5)
holds if and only if $\tilde u=0$, which is equivalent to that $u$
is an affine linear function.  \qed
\enddemo

\demo {Proof of Theorem 0.1} Let $u\ge 0$ be a rational
PL-function which associates to  a toric degeneration on $M$. Then
by Lemma 4.3, we have
 $$\aligned & L(u)> 0,~\text{if}~ u~\text{is not an affine linear function},~\text{or}~\\
 &L(u)=0,~\text{if} ~u~\text{is  an affine linear
 function}.\endaligned$$
The later implies that $u\equiv 0$ since $u\ge 0$, and
consequently, the corresponding  toric degeneration is trivial.
Thus by Proposition 3.1, we have
$$F_{\tilde{\beta}}(\tilde{\alpha})=-\frac{1}{2Vol(P)}L(u)<0,$$
 if the toric degeneration is not trivial, where $F_{\tilde{\beta}}(\tilde{\alpha})$ denotes
 the modified Futaki
invariant on the central fibre associated to a $\Bbb C^*$-action
induced by the toric degeneration. The  theorem is
proved.\qed\enddemo

Next we show that  the $K$-semistability is a necessary condition
for the existence of extremal metrics on a polarized toric
manifold. We need

\proclaim{Lemma 4.4}If $\mu(\phi)$ is bounded from below in $\Cal
M_{G_0}$, then
$$L(u)\geq 0, \forall ~u \in \Cal C_\infty.\tag 4.6$$
\endproclaim

\demo{Proof} Suppose that there exists $f$ in $\Cal C_{\infty}$
such that
$$L(f)<0.$$
Then we choose a function $u$ in $\Cal C$ and consider a sequence
of $u_k=u+kf\in \Cal C_{\infty}$. By the convexity of $\Cal F(u)$,
we have
$$\Cal F(u_k)\leq \Cal F(u)+ kL(f) \rightarrow -\infty$$
as $k \rightarrow \infty$. This is  impossible since $\mu(\phi)$
is bounded from below. Thus (4.6) is true. \qed
\enddemo

\proclaim{Proposition 4.1} Let $M$ be a polarized  toric manifold
which  admits an extremal metric. Then $M$ is $K$-semistable
relative to the $\Bbb C^*$-action induced by $X$ for any toric
degeneration.
\endproclaim

\demo{Proof} Note that the existence of extremal metrics implies
that $\mu(\phi)$ is bounded from below in $\Cal M_{G_0}$.  Then by
Lemma 4.4, (4.6) is true.  Now suppose that the proposition is not
true. Then by Proposition 3.1, there exists a positive PL-function
$u$ such that $L(u)<0$.  By making a small perturbation to $u$,
one will get a smooth function $u'$ in $\Cal C_{\infty}$ such that
$L(u')<0.$ But this is impossible by (4.6). Thus  the proposition
is true. \qed\enddemo

\demo {Proof of Theorem 0.2} By Proposition 2.1, it suffices to
prove that the condition (2.4) holds for any normalized function
$u\in \Cal C_\infty$ with
 $$ \int_{\partial P} u d\sigma=1.\tag 4.7$$
  As in Lemma 4.2, we may assume that $P$ contains the origin $0$ and $u$ is normalized at
  $0$. Moreover,  we have
  $$L(u)\ge 0,$$
 otherwise, one can find a PL-function $u'$ such that
 $$L(u')<0,$$
  which is  impossible according to Lemma 4.3.
  Thus  by  the contradiction, if (2.4) is not true, then
there is a sequence of normalized functions $\{u^{(k)}\}$ in $\Cal
C_\infty$ such that
$$\int_{\partial P} u^{(k)} d\sigma=1\tag 4.8$$
  and
$$L(u^{(k)})\longrightarrow 0, \text{as}~~k\longrightarrow \infty.\tag 4.9$$

By Lemma 4.3 and (4.9), we have
$$L(u^{(k)})\geq\sum_{i=1}^d\int_{P_i}(\frac{n+1}{\lambda_i}-\bar{R}-\theta_X)u^{(k)}dx\longrightarrow
0.$$
 On the other hand, by (4.8), we see that there exists a subsequence (still denoted by $\{u^{(k)}\}$) of
 $\{u^{(k)}\}$, which converges locally
uniformly to a normalized convex function $u_\infty\ge 0$ on $P$.
It follows
$$\aligned &\sum_{i=1}^d \lim_k\int_{P_i} (\frac{n+1}{\lambda_i}-\bar{R}-\theta_X)
u^{(k)}dx\\
&=\sum_{i=1}^d\int_{P_i} (\frac{n+1}{\lambda_i}-\bar{R}-\theta_X)
u_\infty dx=0.\endaligned\tag 4.10$$
 We claim that
  $$u_\infty\equiv 0, ~\text{on}~ P.\tag 4.11$$
    If (4.11) is not true,
 then there  is an open set $U$ of $P$ such that $u_\infty>0$.  Since
each set  $\{\frac{n+1}{\lambda_i}-\bar{R}-\theta_X=0\}\cap P$
lies in a hyperplane in $\Bbb R^n$, we get
 $$\aligned &\sum_{i=1}^d \int_{P_i} (\frac{n+1}{\lambda_i}-\bar{R}-\theta_X)
u_\infty dx\\
&\ge \sum_{i=1}^d \int_{U\cap P_i}
(\frac{n+1}{\lambda_i}-\bar{R}-\theta_X) u_\infty
dx>0,\endaligned$$
 which is a contradiction to  (4.10). Thus (4.11) is true.

  By the assumption conditions (4.8) and (4.9), one sees
 $$\aligned &\lim_k \int_P (\overline R+\theta_X)u^k dx\\
&=\lim_k  \int_{\partial P} u^k d\sigma
 =1.\endaligned$$
On the other hand, similarly to  (4.10), we have
$$\aligned &\lim_k \int_P (\overline R+\theta_X)u^k dx\\
&=\int_P (\overline R+\theta_X)u_\infty dx=0.\endaligned$$
  Thus we get a contradiction from the above. The contradiction
shows that (2.4) is true. The theorem is proved. \qed
\enddemo

\demo{Proof of Corollary 0.1} By (1.9) in Section 1 and Lemma 4.1,
we see
$$ \int_{\partial P} d\sigma= n\sum_{i=1}^d \int_{P_i} \frac {1}{\lambda_i} dx.$$
Then condition (0.2'') in Corollary 0.1 implies condition (0.2)
since $\theta_X=0$. Thus the corollary follows from Theorem
0.2.\qed\enddemo

 For the case of toric Fano surfaces, condition (2.4) can be proved in
another way. Recall that a simple PL-function is a form of
$$u=\max \{0,\sum a_i x_i +c\}$$
for some vector  $(a_i) \in \Bbb R^n$ and number $c \in \Bbb R$.
 We call the
hyperplane $\sum a_i x_i +c=0$  a crease of $u$. The following
result was proved in [D2] as a  special case of Proposition 5.2.2
and 5.3.1 there.

\proclaim{Lemma 4.5 }Let $M$ be a  toric surface  with
$$\bar{R} + \theta_X\ge 0.\tag 4.12$$
  Suppose that the condition (2.4) in
Section 2 doesn't hold. Then either

 1. there exists a positive rational PL-function $u$ such that $L(u) < 0$, or

2.  there exists a simple PL-function $u$ with crease intersecting
the interior of the convex polytope such that $L(u)= 0$.
\endproclaim

 It is known that toric Fano surfaces are classified into five
different types, i.e., $\Bbb CP^2$, $\Bbb CP^1\times \Bbb CP^1$
and $\Bbb CP^2\# l\overline {\Bbb CP^2}~(l=1,2,3)$ and there
exists a K\"ahler-Einstein on each one except $\Bbb CP^2\#
l\overline {\Bbb CP^2}~(l=1,2)$ ( [T1]).  For the last two cases
we will verify the  condition (4.12) in the next section. Thus
combining Lemma 4.3 and 4.5, we see  that condition (2.4) is true
if (4.3) holds.

The condition (4.3) is  also true for toric Fano surfaces (cf.
Section 5). Thus as a corollary of Theorem 0.2, we have

\proclaim{Corollary 4.1} Let $M$ be a toric Fano surface and $G_0$
be a maximal compact subgroup of torus actions group $T$ on $M$.
Then there exists numbers $\delta
>0$ and $C$ such that for any $G_0$-invariant potential functions associated to the K\"ahler class
$2\pi c_1(M)$, it holds
$$\mu(\phi) \geq \delta \inf_{\tau\in {T}}{I(\phi_\tau)} - C_{\delta}.$$
\endproclaim

According to [WZ],  any  toric Fano manifold  admits a
K\"ahler-Einstein metric iff  the Futaki invariant vanishes. Then
by Theorem 0.2, we have

 \proclaim {Corollary 4.2} Let $M$ be a toric Fano manifold with
 the vanishing Futaki invariant. Then there exists  $\delta
>0$ and $C$ such that
$$\mu(\phi) \geq \delta I(\phi) - C_{\delta} \tag 4.13$$
  for any $\phi\in\Lambda_1^{\bot}(M,g_{KE})\cap \Cal M_{G_0},$
 where  $\Lambda_1(M,g_{KE})$ denotes the first non-zero eigenfunctions space of
 the Lapalace operator associated to the $G_0$-invariant K\"ahler-Einstein metric $g_{KE}$ on $M$.
\endproclaim

\demo{Proof} We choose the K\"ahler-Einstein metric $g_{KE}$ on
$M$ as an initial metric in the definition of potential functions
space $\Cal M$ and define  a functional on the automorphisms group
$\text{Aut(M)}$ on $M$ by
$$\Phi(\tau)=(I-J)(\phi_{\tau})$$
  for any $\phi\in\Cal M$,  where $\phi_{\tau}$ is defined as in
(2.8) in Section 2 while $\omega_g$ is replaced by
 $\omega_{KE}$. Then according to [BM], there is a
 $\sigma\in\text{Aut(M)}$ such that
 $$\Phi(\sigma)=\inf_{\tau\in\text{Aut(M)}}\Phi(\tau)$$
 and consequently $\phi_{\sigma}\in\Lambda_1^{\bot}(M,g_{KE})$.
Note that $\phi_{\sigma}$ is invariant under a maximal compact
subgroup of $\text{Aut(M)}$ as same as $\phi$.
 The inverse is also true. In fact, from the proof of uniqueness of
 K\"ahler-Einstein metrics in [BM], one can prove that
$\phi\in\Lambda_1^{\bot}(M,g_{KE})$ iff
 $$ (I-J)(\phi)=\inf_{\tau\in\text{Aut(M)}}\Phi(\tau).$$
  Thus by  Theorem 0.2 and the
inequalities ([Au]),
 $$ \frac{1}{n} J(\phi)\le(I-J)(\phi)\le\frac{n-1}{n}J(\phi),~\forall~\phi\in \Cal M,$$
 we  get
 $$\aligned \mu(\phi)&\ge \delta \inf _{\tau}I(\phi_\tau) - C_{\delta}\\
 &\ge\delta \inf _{\tau}(I-J)(\phi_\tau) - C_{\delta}\\
 &=\delta (I-J)(\phi) - C_{\delta}\\
&\ge \frac{\delta}{2n} I(\phi) - C_{\delta}.\endaligned$$
 Replacing $\frac{\delta}{2n}$ by $\delta$, we obtain (4.13).
\qed\enddemo

\proclaim{Remark 4.1} The existence of such an inequality (4.13)
in Corollary (4.2) was conjectured by Tian for any
K\"ahler-Einstein manifold ([T2]).
\endproclaim

\subheading{\bf 5. $\theta_X$ on $\Bbb CP^2\# \overline{\Bbb
CP^2}$ and $\Bbb CP^2\# 2\overline{\Bbb CP^2}$ }
  \vskip3mm

 In this section, we estimate the quality $\theta_X$ on the
 manifolds $\Bbb CP^2\# 1\overline{ \Bbb CP^2}$ and  $\Bbb CP^2\#$
 \flushpar $2\overline{ \Bbb CP^2}$ and
show that the conditions (4.3) and (4.12) in Section 4 are true
for these two manifolds.

\proclaim{Proposition  5.1}On $\Bbb CP^2\# 2 \overline{\Bbb CP^2}$
or $\Bbb CP^2\# 1 \overline{\Bbb CP^2}$, we have
$$-2 <\theta_{X}<1.$$
\endproclaim

\demo{Proof} First we consider  $\Bbb CP^2\# 2 \overline{\Bbb
CP^2}$. In this case, we choose a
 K\"ahler form
 $$\omega_g =\sqrt{-1}\partial\bar{\partial}\varphi~\in~2\pi c_1(M) $$
 which is given by a convex function,
$$\varphi =
\frac{1}{2}\left[log(1+e^{2y_{1}}+e^{2y_{2}})+ log(1+e^{2y_{1}})+
log(1+e^{2y_{2}})-2y_{1}-2y_{2}\right].$$
 Then $M$ associates  to a pentagon $P_{0}$ in $\Bbb R^{2}$ by using the moment map
$x_{i}=\frac{\partial \varphi}{\partial y_{i}}$. Namely,
$$\cases & \frac{\partial \varphi}{\partial y_{1}}
=\frac{e^{2y_{1}}}{(1+e^{2y_{1}}+e^{2y_{2}})}
+\frac{e^{2y_{1}}}{(1+e^{2y_{1}})}-1\\
& \frac{\partial \varphi}{\partial y_{2}}
=\frac{e^{2y_{2}}}{(1+e^{2y_{1}}+e^{2y_{2}})}
+\frac{e^{2y_{2}}}{(1+e^{2y_{2}})}-1.\endcases $$
 Thus $P_{0}$
is decided around by five edges:
$$x_{1}=1,\ x_{2}=1,\ x_{1}=-1,\ x_{2}=-1,\ x_{1}+x_{2}=1.$$

 For any holomorphic vector field $v$ on $M$, we normalize
the potential function $\theta_{v}$ of $v$ by
$$\int_M \theta_{v}e^{h_g}\omega_g^n=0,$$
 where $h_g$ is a potential function of the Ricci form of
 $\omega_g$.   Then $\theta_{v}$ satisfies (see [TZ]),
$$\theta_{v}=-\triangle \theta_{v}- v(h_g).$$
So the Futaki invariant can be computed by
$$F(v)=\int_M v(h_g)\omega^n=-\int_M \theta_{v}\omega^n=
-(2\pi)^2\int_{P_0}\theta_{v}dx.$$
 In particular,
$$\aligned
F_{1}=F(\frac{\partial}{\partial y_1})
& = -(2\pi)^2\int_{P_{0}} x_{1} dx_{1}dx_{2} \\
& = -\frac{(2\pi)^2}{3},
\endaligned\tag 5.1$$
and
 $$F_{2}=F(\frac{\partial}{\partial y_2})=-\frac{(2\pi)^2}{3},\tag 5.2$$
  where $(y_1,y_2)$ is   an affine coordinates system.

  By using a translation
$x'_{i}=x_{i}+\frac{1}{3vol(P)}$, we get a new pentagon $P_{1}$ so
that
$$ \int_{P_{1}} x'_{i} dx'_{1}dx'_{2} = 0,~i=1,2.$$
Then  $\theta_X$ associated to the extremal holomorphic vector
field $X$ is a from of
$$\theta_{X} =\theta_{1}x'_{1}+\theta_{2}x'_{2},$$
 where $\theta_1$ and $\theta_2$ are two constants.

Since the Futaki invariant in the system $(x'_1,x'_2)$ is computed
  by ([FM]),
$$\aligned F(v)&=\int_M v(h_g)\omega_g^n=-\int_M \theta_{v}\triangle
h_g\omega_g^n =-\int_M
\theta_{v}\theta_{X}\omega_g^n\\
 &=-(2\pi)^2\int_{P_1}\theta_{v}\theta_{X}dy,\endaligned $$
 we get from (5.1) and (5.2),
$$\cases & \left( \int_{P_{1}} {x'}^{2}_{1} dx'_{1}dx'_{2} \right) \theta_1 +
\left( \int_{P_{1}} x'_{1}x'_{2} dx'_{1}dx'_{2} \right) \theta_2
=  \frac{1}{3}\\
& \left( \int_{P_{1}} x'_{1}x'_{2} dx'_{1}dx'_{2} \right)
\theta_{1} + \left( \int_{P_{1}} {x'}_{2}^{2} dx'_{1}dx'_{2}
\right) \theta_{2} = \frac{1}{3}.\endcases $$
 A simple computation shows
$$\theta_{1}=\theta_{2}=-\frac{168}{409}.$$
  Thus $-2<\theta_{X} <1$.

 In the  case of $\Bbb CP^2\# 1 \overline{\Bbb CP^2}$,  we choose a
 K\"ahler form
 $$\omega_g =\sqrt{-1}\partial\bar{\partial}\varphi~\in~2\pi c_1(M) $$
 which is given by a convex function,
$$\varphi =
\frac{1}{2}\left[2log(1+e^{2y_{1}}+e^{2y_{2}})+
log(e^{2y_{1}}+e^{2y_{2}})-2y_{1}-2y_{2}\right].$$
 So $M$
associates to a quadrilateral $P_{0}$ in $\Bbb R^{2}$ decided
around by four edges:
$$x_{1}+x_{2}=-1,\ x_{1}=-1,\ x_{2}=-1,\ x_{1}+x_{2}=1.$$
It can be computed in the similar way that $$\theta_{X}
=\frac{5}{29}x'_{1}+\frac{5}{29}x'_{2}$$ under the coordinates
$\{x'_{1},x'_{2}\}$ satisfying
$$ \int_{P_{1}} x'_{i} dx'_{1}dx'_{2} = 0,~i=1,2.$$
Thus $-2<\theta_{X}<1$. \qed
\enddemo

\subheading{\bf 6. K\"ahler classes satisfying  (0.2'') on $\Bbb
CP^2\# 3\overline{\Bbb CP^2}$}

 \vskip3mm
In this section, we construct a kind  of K\"ahler classes  on  the
toric manifold $\Bbb CP^2\# 3\overline{\Bbb CP^2}$, which satisfy
the condition (0.2'')  in Introduction. It is well known that
$\Bbb CP^2\# 3\overline{\Bbb CP^2}$  with an anti-cannonical line
bundle is corresponding to an integral polytope in $\Bbb R^2$,
 $$\aligned P_0=\{( x_1,x_2)\in \Bbb R^2|~ &x_1\leq 1,
 -x_2\leq 1,  -x_1-x_2\leq 1, \\
& -x_1\leq 1, x_2\leq 1, x_1+x_2\leq 1\}. \endaligned$$
 Let
$$l_1: x_1 = 1; $$
$$ l_2: -x_2 = 1; $$
$$ l_3: -x_1-x_2=1;$$
$$l_4: -x_1 =1; $$
$$l_5: x_2 = 1; $$
$$l_6: x_1+x_2=1; $$
 be six edges of $P_0$,  and $D_i$ be six divisors  on $\Bbb CP^2\# 3\overline{\Bbb CP^2}$ corresponding to
$l_i$. Then the cohomology  class of a K\"ahler class is given by
$$\sum_{i=1}^6 \lambda_i [D_i], $$
  where  $\lambda_i$ are all positive numbers. So the K\"ahler class is corresponding to a convex
polytope
$$\aligned P=\{ (x_1,x_2)\in \Bbb R^2|~ &x_1\leq \lambda_1,
 -x_2\leq \lambda_2,
  -x_1-x_2\leq \lambda_3,\\
& -x_1\leq \lambda_4,  x_2\leq \lambda_5,  x_1+x_2\leq
\lambda_6\}.\endaligned$$
 Moreover  the K\"ahler class is integral when  $\lambda_i$ are all positive
integers.

\proclaim {Lemma 6.1} Let $\lambda$ and $\mu$ be two positive
numbers with
 $\frac{\lambda}{2}<\mu<2\lambda$. Let
 $$\lambda_1=\lambda_3=\lambda_5=\lambda,~
\lambda_2=\lambda_4=\lambda_6=\mu.$$
 Then the Futaki invariant $F(.)$ vanishes on  the K\"ahler class  associated to $P$.\endproclaim

\demo{Proof} By a direct computation, one sees
$$\int_P x_idx=0~\text{and}~\int_{\partial P} x_idx=0,~i=1,2.$$
 Then
$$F(\frac{\partial}{\partial y_i})= 2^n n!(2\pi)^n(\int_{\partial P} x_i-\int_P
\overline R x_i dx) =0, ~i=1,2.$$
 Thus the Futaki invariant vanishes.\qed\enddemo

\proclaim {Lemma 6.2} Let  $\lambda$ and $\mu$, and $\lambda_i,
~i=1,...,6$ be numbers as in Lemma 6.1.   Suppose
$$\frac{\lambda}{1+\frac{\sqrt{10}}{5}} \le \mu \le (1+\frac{\sqrt{10}}{5})\lambda. \tag 6.1$$
 Then condition (0.2'')  is satisfied  on
 the K\"ahler class  associated to
$P$ .
\endproclaim

\demo{Proof} By (1.9) in Section 1, one sees
 $$\bar{R}=\frac{2(\mu+\lambda)}{4\lambda\mu-\mu^2-\lambda^2}.$$
 Then condition (0.2'') is equivalent to
$$\frac{2(\mu+\lambda)}{4\lambda \mu-\mu^2-\lambda^2}\le \frac{3}{\text{max}\{\lambda,\mu\}}.$$
The late is equivalent to (6.1).\qed\enddemo

Combining  Lemma 6.1 and 6.2 to Theorem 0.2, we have

\proclaim{Proposition 6.1} Let  $\lambda$ and $\mu$ be two
positive numbers which satisfy (6.1). Let
 $$\lambda_1=\lambda_3=\lambda_5=\lambda,~
\lambda_2=\lambda_4=\lambda_6=\mu.$$
 Then on the K\"ahler class with cohomology class
$$ \lambda ( [D_1] + [D_3]+ [D_5])+
\mu ( [D_2] + [D_4]+ [D_6]),$$
 the $K$-energy $\mu(\phi)$ is proper associated to subgroup $T$ in the space of
$G_0$-invariant K\"ahler metrics  on $\Bbb CP^2\# 3\overline{\Bbb
CP^2}$ associated to $P$.
\endproclaim

It is clear that for a pair $(\lambda,\mu)=(2,3)$ or
$(\lambda,\mu)=(3,2)$, (6.1) is satisfied. Then  by Proposition
6.1, $K$-energy $\mu(\phi)$ is proper  on the K\"ahler class  on
$\Bbb CP^2\# 3\overline{\Bbb CP^2}$ associated to the integral
polytope $P$ with
$$\lambda_1=\lambda_3=\lambda_5=2,~ \lambda_2=\lambda_4=\lambda_6=3,$$
 or
$$\lambda_1=\lambda_3=\lambda_5=3,~ \lambda_2=\lambda_4=\lambda_6=2.$$
It is interesting to study whether there exists  a K\"ahler metric
with constant scalar curvature on  $\Bbb CP^2\# 3\overline{\Bbb
CP^2}$ in these two K\"ahler classes.

\enddocument